\documentclass[12pt,leqno]{amsart}
\usepackage{latexsym,amssymb,amsmath,multicol, rotating,lscape}
\usepackage{letltxmacro}

\usepackage{comment}
\usepackage{mwe}
\usepackage{amsmath}

\usepackage{tikz}
\usetikzlibrary{arrows.meta, positioning}
\usetikzlibrary{decorations.markings}
\tikzset{
    process/.style={
        text width=4cm, draw,
        minimum height=1.6cm,
        text centered,
        },
    description/.style={
        text centered,
        text width=10cm,
    },
    myarrow/.style={
        postaction={
            decorate, decoration={
                markings,mark=at position #1 with {\arrow{Stealth};
                }
            }
        }
    },
}
\tikzset{
    process1/.style={
        text width=3cm,
        minimum height=1cm,
        text centered,
        },
    description/.style={
        text centered,
        text width=10cm,
    },
    myarrow/.style={
        postaction={
            decorate, decoration={
                markings,mark=at position #1 with {\arrow{Stealth};
                }
            }
        }
    },
}
\usepackage{setspace}
\usepackage{etoolbox}
\AtBeginEnvironment{tikzpicture}{\singlespacing}
 2

\newtheorem{thm}{Theorem}[section]

\newtheorem{conjecture}[thm]{Conjecture}

\newtheorem{cor}[thm]{Corollary}
\newcommand{\thmref}[1]{Theorem~\ref{#1}}

\theoremstyle{remark}

\begin{document}

\title[Elliptic curve over totally real number fields : A Survey]
{Elliptic curve over totally real fields: A Survey}

\author{Bidisha Roy}
\address[Bidisha Roy]{Scuola Normale di Pisa, Piazza dei Cavlieri, 56126, Pisa, Italy}
\email{brroy123456@gmail.com}
\author{ Lalit Vaishya}
\address[Lalit Vaishya]{ The Institute of Mathematical Sciences, CIT Campus, Taramani,Chennai-600 113
Tamilnadu, India}
\email{lalitvaishya@gmail.com, lalitv@imsc.res.in}

\begin{abstract}
    In this survey article, we summarise the known results towards the conjecture: elliptic curves over totally real number fields are modular. For  understanding these recent results in the literature, we present some necessary background  along with certain applications. 
\end{abstract}

\maketitle

\tableofcontents

\section{Introduction}
  Let us start this survey-article by recalling the following  conjecture which is a generalized version of Shimura - Taniyama Conjecture (1960) (currently, Theorem). \\
  \noindent\textbf{ Conjecture:} ~ {\it Every elliptic curve defined over  a number field $K$ is  modular.}   \\
To be precise,  an elliptic curve $E$ over a totally real number field $ K$ of conductor $N$ is called {\it modular} when there exists a Hilbert newform $\mathfrak{f}$ of level $N$, parallel weight $2$ and with rational Hecke
eigenvalues, such that their corresponding systems of compatible $p$-adic Galois representations are isomorphic,$$ \rho_{E,p^{\infty}} \cong \rho_{\mathfrak{f},p},$$
where the left  side is $p$-adic Galois representation attached to $E$ arises from the action of $Gal(\overline{K}/K)$ on the Tate module $T_p(E)$, and the right side is the Galois representation associated to a Hilber modular form $\mathfrak{f}$  (defined
by Taylor \cite{Taylor1989}). We  devote this article to understand all the terms mentioned in the above lines. 

\smallskip

First, it is natural to understand an elliptic curve $E$ over $\mathbb{Q}$. In that case,  one needs to look for   a cuspidal Hecke eigenform $ f$ on some congruence subgroup $ \Gamma_0(N) \subset SL_2(\mathbb{Z})$ such that the Galois representation on the $p$-adic Tate-module of $ E$ and  the $ p$-adic representation attached to $f$ (see \cite{DSbook}) are isomorphic.  In \cite{Ribet}, Ribet’s result (on $\epsilon$-conjecture of Serre \cite{Serre})  tells us that in order to deduce Fermat’s Last Theorem, it is enough to prove the above conjecture for semistable elliptic curves over $\mathbb{Q}$. Wiles  approached to study modularity conjecture of  elliptic curves over $\mathbb{Q}$ via Galois representations. With upgradation in Talyor-Wiles \cite{Taylor-Wiles}, they could prove that semistable  elliptic curves over $\mathbb{Q}$  are modular, which leads to proving Fermat's Last Theorem.   The full conjecture for $ K = \mathbb{Q}$ is known due to Wiles \cite{Wiles}, Breuil, Conrad, Diomand and Taylor \cite{BCDT}.  It is known as {\it modularity theorem over $\mathbb{Q}$} in the literature. 

\smallskip

Though it is natural to understand the generalization of the conjecture to general number fields (in other words, to show all elliptic curves over $K$ have compatible systems of $p$-adic representations of $Gal(\overline{K}/K)$ corresponding to a cuspidal automorphic representation $\pi$ of $ GL_2(\mathbb{A}_K)$)  but  the existence of Galois representations associated to $\pi$ is not known.  For a totally real field, the Galois representations have been achieved in the literature. With the help of these developments, the modularity of elliptic curves over some families of  number fields are present now in literature.  The first modularity theorem, the Fermat equation over $ \mathbb{Q}(\sqrt{2})$, has been obtained by Jarvis and Meekin in 2004 \cite{JM2004}. This result relied on  Fujiwara's \cite{Fujiwara} unpublished result. The first unconditional result, modularity of semistable elliptic curves over $ \mathbb{Q}(\sqrt{2})$ and $ \mathbb{Q}(\sqrt{17})$, was obtained by  Jarvis and Manoharmayum   (in \cite{JM2003, JM2008}). After that, a major achievement was obtained  by Freitas, Hung and Siksek \cite{FHS2015} when they established the modularity of elliptic curves over all real quadratic fields. Later, Derickx, Najman and Siksek \cite{DNS2020} extended it to elliptic curves over all totally real cubic fields. More recently,  Box \cite{Box2021} could prove that every elliptic curve defined over a totally real quartic number field not containing $\sqrt{5}$ is modular.

\smallskip

In the local set-up, there are various other works which establish modularity under local assumptions on the elliptic curve $E $ over number field $ K$. In a particular direction, Dieulefait-Freitas  \cite{DF2013} considered modularity of the curve $x^{13} + y^{13} = C z^p$, for $p > 13649$ and for an infinite family of values for $C$. These kinds of problems over local fields have been considered by other number theorists also (for instance, \cite{Bruin2000, CS2009, thorne2022}).

\smallskip

In this article, we  survey the known results in the direction of elliptic curves over totally real fields. More concretely, we mainly discuss about elliptic curves over totally real quadratic and cubic fields because number theorists  first achieved these two general results in this wide open-direction. Our goal is to understand the  hurdles in obtaining the modularity result over number fields in  comparison  to the rational case and to get an overview of these obstacles and produce a comparative study. 

\smallskip

The arrangements of this article is as follows. In section 2, we provide the necessary background to  understand the known results in the literature. In this section, we also briefly recall Hilbert modular forms and its associated Galois representation.  In Section 3, we sketch the proofs of the main results associated with the modularity of an elliptic curve over certain number fields. In Section 4, we conclude by briefing some recent applications of the modularity results.

 \section{Background}

Let $K$ be a number  field. We consider an elliptic curve $E$ defined over $K$ in a Weiestrauss form given by;
$$E: \quad y^2 = x^3 +Ax+B, ~{\rm where} ~ A, B \in K ~ { \rm ~with~ discriminant ~}\Delta= -16(4A^{3}+27B^{2}).$$ 
In particular, for $K=\mathbb{Q}$, an elliptic curve $E$ over $\mathbb{Q}$ is said to be modular if it has a finite covering by a modular curve  $X_{0}(N)$. Such an elliptic curve has the property that its  associated Hasse-Weil zeta $L$-function has an analytic continuation and satisfies a nice functional equation of the standard type. If an elliptic $E$ curve over $\mathbb{Q}$ with a given $j$-invariant is modular, it is easy to see that all elliptic curves with the same $j$-invariant are modular.   The key development in the proof of modularity theorem for semistable  elliptic curves over $\mathbb{Q}$  shows  a  link between Galois representations and modular forms on the one hand and the interpretation of special values of $L$-functions on the other (see for detail Teorem \ref{modularity over Q} and thereafter.) Wiles mentioned that  this method seems well suited to establish the modularity theorem for all elliptic curves over $\mathbb{Q}$ and it can be  generalized to other totally real number fields. 

\smallskip
 
 Here, we produce brief details of Wiles's approach towards modularity theorem for elliptic curve over $\mathbb{Q}$. Let $E$ be an elliptic curve defined over $\mathbb{Q}$ with discriminant $\Delta$. For a rational prime  $p$, $\mathbb{Z}/p \mathbb{Z}$ is denoted as a finite field of characteristic $p$. Let  $a_{p} := p+1 - \#E(\mathbb{Z}/p \mathbb{Z})$, where
$$ E(\mathbb{Z}/p \mathbb{Z}) := \{ (x,y) \in \left(\mathbb{Z}/p \mathbb{Z}\right)^{2} : y^2 \equiv x^3 +Ax+B (\rm mod ~p) ~{\rm has ~ a ~solution  } \}.$$ 
With the above notations, the  Hasse-Weil zeta $L$-function is defined as follows:
\begin{equation*}
\begin{split}
L(E/\mathbb{Q},s) = \displaystyle{\prod_{p \mid \Delta }\left(1- \frac{a_{p}}{p^{s}} \right)^{-1}
\prod_{p \nmid \Delta }\left(1- \frac{a_{p}}{p^{s}} +\frac{1}{p^{2s-1}} \right)^{-1} }. 
\end{split}
\end{equation*} 
The Hasse-Weil zeta $L$-function converges absolutely for $\Re(s)>3/2$. Let $f$ be a classical normalised  Hecke eigenform (a modular cusp fom) of weight $2$ and level $N$ with the Hecke eigenvalues $a_{f}(n)$. Associated to $f$, we define Hecke $L$-function given by
$$L(f,s)= \displaystyle{\sum_{n=1}^{\infty}\frac{a_{f}(n)}{n^{s}} }, ~ \qquad  \Re(s)>3/2.$$
Now, we are ready to state the following modularity theorem over $\mathbb{Q}$.

\begin{thm}
{\rm \textbf{(Modularity theorem over $\mathbb{Q}$):} } \label{modularity over Q}
Let $E$ be an elliptic curve defined over  $\mathbb{Q}$. Then, there exist a  classical Hecke eigenform $f$ of weight $2$ and some level $N$ such that 
$$
L(E/\mathbb{Q},s) = L(f,s) ~ \text{for all}~ s \in \mathbb{C}.
$$
\end{thm}

As we mentioned earlier, this was proved  using Galois representation theoretic approach. In a different set-up (a representation-theoretic approach), the modularity theorem over $\mathbb{Q}$ states that  every elliptic curve defined over $\mathbb{Q}$ is modular if there exists a Hecke eigenform of weight $2$ and level $N$ such that their corresponding system of compatible $p$-adic Galois representations are isomorphic. Here, the $p$-adic Galois representation for $E$ over $\mathbb{Q}$ arises from the action of $Gal( \overline{\mathbb{Q}}/\mathbb{Q})$ on Tate module $T_{p}(E) = {\varprojlim}~ E[p^{n}]$, a projective limit of the set of  $p^{n}$-division points of an elliptic curve $E$ over $\mathbb{Q}$. Now, we set up the preliminaries to state the precise representations theoretic statement of modularity theorem  over $\mathbb{Q}$ for an elliptic curve $E$ over $\mathbb{Q}$. Some of these preliminaries are available over any number field $K$ and we state those over $K$.

\smallskip

\subsection{mod-n Galois representation}: Let $E(K)$ denote the set of $K$-rational points containing a point at infinity (say $\mathcal{O}$). It forms an abelian group under some formal group law. A celebrated theorem of Mordell and  Weil tells that  $E(K)$ is a finitely generated abelian group over the number field $K$. More precisely, 
$$
E(K) \cong \mathbb{Z}^{r_{alg}} \oplus {\rm ~torsion ~ part}
$$
where $r_{alg}$ is known as {\it algebraic rank of an elliptic curve $E$ over $K$}.
 For  a positive integer $n$, let $E[n]$ denote the $n$-torsion subgroup of $ E(\overline{K})$, where $ \overline{K}$ is an algebraic closure of $K$. It is well-known fact that $E[n]$ is a free module of rank $2$ over $\mathbb{Z}/n \mathbb{Z}$. For a given positive integer $n$, the action of $Gal( \overline{K}/K)$ on the torsion points of $ E[n]$ is well known in the literature which is actually the action of $Gal( \overline{K}/K)$ on the coordinates of the points on $E$. It gives rise to a {\it mod $n$ Galois representation} (attached to $E$) given by 
 
$$
\overline{\rho}_{E, n}: Gal( \overline{K}/K) \longrightarrow Aut (E[n])
$$
where $ Aut (E[n])$ denotes the group of automorphism of abelian group $E[n]$. Since $E[n]$ is isomorphic to the group $\left(\mathbb{Z}/n \mathbb{Z}\right)^{2}$. This leads $ Aut (E[n]) \cong GL_{2}\left(\mathbb{Z}/n \mathbb{Z}\right)$. Each choice of basis determines such an isomorphism and various isomorphism obtained in this process differ only by inner automorphism. Therefore, each element of $ Aut (E[n])$ has a well defined trace and determinant in $\mathbb{Z}/n \mathbb{Z}$. Moreover, fixing a choice of basis of  $\left(\mathbb{Z}/n \mathbb{Z}\right)^{2}$, a representation can viewed as taking the values from the set $GL_{2}\left(\mathbb{Z}/n \mathbb{Z}\right)$ like;  $$
\overline{\rho}_{E, n}: Gal( \overline{K}/K) \longrightarrow GL_{2}\left(\mathbb{Z}/n \mathbb{Z}\right).
$$

The kernal of this representation corresponds to a finite Galois extension $K_{n}$  of $K$ in  $\overline{K}$ via Galois theory. Concretely,  let $G_{n} := Gal (K(E[n])/K)$. For each prime $p$, not divining $n$ and ${\rm conductor(E)}$, a familiar construction produces a Frobenius  element $\sigma_{p}$ in $G_{n}$ which is well defined up to conjugation. More precisely,  $\overline{\rho}_{E,n} (Frob_{p}) = \sigma_{p}$. On the other hand, we have  ${\rm trace}(\sigma_{p}) \equiv a_{p} (\rm mod~ n), $ where $a_{p} = p+1 - \#E\left(\mathbb{Z}/p \mathbb{Z}\right)$. Thus, the representation $\overline{\rho}_{E, n}$ captures the information about the number $a_{p}$. In the later part of this article, we shall witness the importance of $ \overline{\rho}_{E,p}$ in more  elaborate way. Now we connect it with another representation as follows.

\smallskip

\subsection{$\ell$-adic Galois representation over $\mathbb{Q}$}: Let $\ell$ be  a rational prime. We  consider the family of groups $E[\ell^{n}]$,  $\ell^{n}$-division points of elliptic curves $E$  for $n =1,2,3, \ldots$. We define the Tate module $T_{\ell}(E)$ as a projective limit of the set of  $\ell^{n}$-division points of an elliptic curve $E$ over $\mathbb{Q}$, i.e., $T_{\ell}(E) = {\varprojlim}~ E[\ell^{n}]$. Moreover, we have $Aut (T_{\ell}(E)) \cong  {\varprojlim} \left( Aut( E[\ell^{n}]) \right)  \cong GL_{2}\left({\mathbb{Z}}_{\ell}\right) $.  For a given rational prime $\ell$,  the resulting sequence of continuous Galois representations 
$
\overline{\rho}_{E, \ell^{n}}: Gal( \overline{\mathbb{Q}}/\mathbb{Q}) \longrightarrow
GL_{2}\left(\mathbb{Z}/\ell^{n} \mathbb{Z}\right) 
$
can be viewed as a single representation 
$$
\rho_{E, \ell^{\infty}}: Gal( \overline{\mathbb{Q}}/\mathbb{Q}) \longrightarrow
GL_{2}\left({\mathbb{Z}}_{\ell}\right) 
$$
where ${\mathbb{Z}}_{\ell}$ is the ring of $\ell$-adic integers. The representation $\rho_{E, \ell^{\infty}}$ is known as a {\it $\ell$-adic Galois representation} attached  to elliptic curve $E$ over $\mathbb{Q}$. Moreover, ${\rm trace}_{\rho_{E, \ell^{\infty}}}({\rm Frob}_{p})$ coincides with the rational integer $a_{p}$ for all prime $p \nmid \ell$ and coprime with the conductor of $E$.

\smallskip

\subsection{Modular representation associated to a Hecke eigenform}:  Let $k\ge 2$  and $\chi$ be a Dirichlet character modulo $N$, for integers $ N \geq 1$. Also let $S_{k}(\Gamma_{0}(N), \chi)$ be the space of cusp forms of weight $k$ for $\Gamma_0(N)$ with character $\chi$. A cusp form $f \in S_{k}(\Gamma_{0}(N), \chi)$ is said to be an eigenform if it a common eigenfunction for all the Hecke operators ($T_{n};n\in \mathbb{N}$) and Atkin-Lehner operators ($W$-operators) \cite[Chapter 5]{DSbook}. Let $f \in S_{k}(\Gamma_{0}(N), \chi) $ be an eigenform. Then, for each $ n\in \mathbb{N}$, there exist an algebraic integer $a_{f}(n)$ such that 
 $T_{n} f = a_{f}(n) f$. Let ${\mathbb{Q}}_{f}$  be the number field generated over $\mathbb{Q}$ by adjoining the elements of the set $\{ a_{f}(n): n\in \mathbb{N}  \}$ and  $\mathcal{O}_{f}$ be its ring of integers.  For a prime $\lambda$ of $\mathcal{O}_{f}$, let $\mathcal{O}_{f, \lambda}$ be the completion of $\mathcal{O}_{f}$ at $\lambda$.  In this set-up, due to  Eichler and Shimura \cite{GS71} (for $k = 2$) and Deligne \cite{Del} (for $k > 2$), we have the following  famous theorem which  connects absolute Galois group over $\mathbb{Q}$ with certain modular representation. 
 \begin{thm}
Let $p$  be a rational prime. Then, for each prime $ \lambda \in \mathcal{O}_{f}$ with $\lambda \mid p$,  there exist  a continuous $\lambda$-adic irreducible Galois  representation 
$$
\rho_{f, \lambda}: Gal( \overline{\mathbb{Q}}/\mathbb{Q}) \longrightarrow GL_{2}(\mathcal{O}_{f, \lambda})  
$$
associated to a normalised Hecke eigenform $f \in S_{k}(\Gamma_{0}(N), \chi)$ such that for all prime $q \nmid Np$,  $\rho_{f, \lambda}$ is unramified at $q$, $
{\rm trace}_{\rho_{f, \lambda}}({\rm Frob}_{q}) = a_{f}(q)$ and $  {\rm det}_{\rho_{f, \lambda}}({\rm Frob}_{q}) = \chi(q) q^{k-1}.
$
\end{thm} 
The analogous result when $k = 1$ is a theorem of Serre and Deligne \cite{SDel} and it is more naturally stated in terms of complex representation.  Such representation $\rho_{f, \lambda}$ is known {\it modular representation associated to a Hecke eigenform $f$}. The modular representation  $\rho_{f, \lambda}$ is odd in the sense that $det({\rho_{f, \lambda}})$ of complex conjugation is $-1$. Moreover, the  modular representation  $\rho_{f, \lambda}$ is potentially-semi stable  at $\ell$ in the sense of Fontaine--Mazur \cite{Font-Maz}.

 \smallskip

  \subsection{ Overview of the modularity theorem over $\mathbb{Q}$}:
Suppose that 
$$ \rho : Gal (\overline{\mathbb{Q}}/ \mathbb{Q}) \rightarrow GL_2(\overline{\mathbb{Q}}_{\lambda})$$ is a continuous representatin. If there exist two positive integers $ k$, $N$  such that there exists an eigenform $ f \in S_k (\Gamma_0(N))$ and  a place $ \lambda$, above a rational prime $p$, then 
$$ \rho \sim \rho_{f, \lambda}$$ holds and we call that {\it $\rho$ is modular.}
Moreover, if a  continuous representation $$ \overline{\rho}:  Gal (\overline{\mathbb{Q}}/ \mathbb{Q}) \rightarrow GL_2(\mathbb{Q}_{p}) $$ lifts to a modular representation, then we say that {\it $ \overline{\rho}$ is modular}.

Suppose $E$ is an elliptic curve defined over the field $\mathbb{Q}$ of rational numbers, and $p$
is an odd prime. The $p$-power torsion points $E[p^\infty] \subset  E(Q)$ form a group isomorphic
to $\left( \mathbb{Q}_p / \mathbb{Z}_p \right)^2$, so the action of the Galois group $Gal(\overline{Q}/Q)$ on $ E[p^\infty]$ defines a
representation $$
\rho_{E, p^{\infty}}: Gal( \overline{\mathbb{Q}}/\mathbb{Q}) \longrightarrow
GL_{2}\left({\mathbb{Z}}_{p}\right) 
$$ as mentioned earlier. Then we have the following theorem due  to  Wiles \cite{Wiles} (for semi-stable elliptic curves over $\mathbb{Q}$ ) and  Breuil, Conrad, Diomand and Taylor \cite{BCDT} ( for any elliptic curves over $\mathbb{Q}$).

\begin{thm}\label{ModularQ} \cite{BCDT}
Let  $E$  be an elliptic curve over $\mathbb{Q}$. Then there exist a Hecke eigenform $f$ of weight $2$ and level $N_{E}$ such that  for all rational prime $p$,  
$$\rho_{E, p^{\infty}} \sim  \rho_{f, \lambda}   \quad {\rm for ~~ some ~~} \quad   \lambda \mid p. $$
Such elliptic curves are known as modular elliptic curves.
\end{thm}

 In this literature, mathematicians first concentrated 
 on semistable elliptic curves due to the link of Fermat's last theorem. So we first focus on this family of curves. Let $E$  be a semistable elliptic curve over $\mathbb{Q}$. Recall from earlier discussion  that  $\rho_{E, p^{\infty}}$  is the representation 
 of $ Gal( \overline{\mathbb{Q}}/\mathbb{Q})$ on the $p$-adic Tate module $T_{p}(E)$ and  $\overline{\rho}_{E, p}$ is the represention of  $ Gal( \overline{\mathbb{Q}}/\mathbb{Q})$ on  the  group of the $p$-division points of $E( \overline{\mathbb{Q}})$. Also consider  $N_{E}$ as the conductor of $E$.   At first one needs to consider $\overline{\rho}_{E,3}$.  The choice of $3$ is critical because a crucial theorem of Langlands-Tunnell \cite{Langlands, Tunnell} shows that if   $\overline{\rho}_{E,3}$ is irreducible then it is also modular. Wiles \cite{Wiles} proved that if  $\overline{\rho}_{E,3}$ is modular then $E$ is modular elliptic curve over $\mathbb{Q}$. 

 \smallskip
 
 Then,  one needs to proceed by showing that under the hypothesis that   $\overline{\rho}_{E,3}$ is semistable at $3$, together with some milder restrictions on the ramification of  $\overline{\rho}_{E,3}$ at the other primes, every suitable lifting of  $\overline{\rho}_{E,3}$ is modular.   To understand the proof of the semistable case, following Taylor-Wiles \cite{Taylor-Wiles},  it is  sufficient to show the modularity of $E$ because it is known that $E$ is modular if and only if the associated $3$-adic representation is modular. 

 \smallskip
 
 The main step to prove \thmref{ModularQ}, for any elliptic curve over $\mathbb{Q}$, is to extend  the methods of Wiles \cite{Wiles} and Taylor-Wiles \cite{Taylor-Wiles} which can be divided into following three cases.

\begin{enumerate}
 \item If $\overline{\rho}_{E, 5}|_{Gal( \overline{\mathbb{Q}}/\mathbb{Q}(\sqrt{5}))}$ is irreducible, then $\rho_{E, 5^{\infty}}$ is modular. 
 \item If $\overline{\rho}_{E,5}|_{Gal( \overline{\mathbb{Q}}/\mathbb{Q}(\sqrt{5}))}$ is reducible but  $\overline{\rho}_{E, 3}|_{Gal( \overline{\mathbb{Q}}/\mathbb{Q}(\sqrt{-3}))}$ is absolutely irreducible, then  $\rho_{E, 3^{\infty}}$ is modular.
 
 \item If $\overline{\rho}_{E,5}|_{Gal( \overline{\mathbb{Q}}/\mathbb{Q}(\sqrt{5}))}$ is reducible but  $\overline{\rho}_{E, 3}|_{Gal( \overline{\mathbb{Q}}/\mathbb{Q}(\sqrt{3}))}$ is absolutely reducible. Then $E$ is isogenous to a modular  elliptic curve with the following $j$-invariant $0, (11/2)^{3}$, or $ -5 (29)^{3}/2^{5}$. So from the examples of the modular elliptic curve of the low conductor, it is modular.  
\end{enumerate} 
 In case (1), irreducibility of $\overline{\rho}_{E, 5}|_{Gal( \overline{\mathbb{Q}}/\mathbb{Q}(\sqrt{5}))}$ implies that $\overline{\rho}_{E, 5}$ is modular and then by showing   modularity of $\overline{\rho}_{E, p}$, one can achive the modularity of $\rho_{E, p^{\infty}}$ with $p=5$. This is known as {\it $3-5$ modularity lifting theorem}. We will see this step in more detail in the following section.  In case (2), It is solved due to the theorem of Langlands and Tunnell \cite{Langlands, Tunnell}. In fact, in both cases $E$ obtain semi-stable reduction over a tame extension of $\mathbb{Q}_{p}$, thus one deduces  the modularity of $\rho_{E, p}$ from the modularity of  $\overline{\rho}_{E, p}$. In the third case, one needs to analyse the rational points on some modular curves of small levels. Thus, one has the modularity of elliptic curve $E$ over $\mathbb{Q}$.

 
 \subsection{Hilbert modular forms and its Galois Representations}
In this section, we produce  a brief overview of the basic theory of Hilbert modular forms and its associated Galois representation.  This is required to explain one main term in the aforementioned  Conjecture mentioned at the beginning of the introduction of this survey.

 \subsubsection{Classical Hilbert Modular forms}
Let $K$ be a totally real number field of degree $n$ and let $\mathcal{O}_K $ be the  ring of integers.  Let  $ \mathbb{H} $ denote the Poincar\'e upper half plane and  $ \sigma_1,\sigma_2,\ldots,\sigma_n $ are the real embeddings of $ K$.  There is a natural action of $ GL_2^+(\mathbb{R})^n $ on $ \mathbb{H}^n $ by component-wise M\"obius transformation.  Given an integral ideal $ \mathfrak{m} $ and a fractional ideal $ \mathfrak{f} $ in $ K $, we define a congruence subgroup
\begin{equation*}\label{Equation "Congruence subgroup definition"}
	\Gamma(\mathfrak{f},\mathfrak{m}) := \left\{ 
	\begin{pmatrix}
		a & b\\
		c & d \\
	\end{pmatrix}
	\in GL_{2}(K) :
	a\in\mathcal{O}_K, b \in \mathfrak{f}^{-1}, c \in \mathfrak{m f}, d\in \mathcal{O}_K, ad - bc \in {\mathcal{O}_K}^{*}
	\right \}.
\end{equation*}
Consider $ {\mathcal{O}_{K_{+}}}^{*} $ as the group of totally positive units of $ \mathcal{O}_{K} $, and $ \chi_0 $ and $ \chi_1 $ are finite order characters of $ (\mathcal{O}_F/\mathfrak{m})^{*} \times (\mathcal{O}_F/\mathfrak{m})^{*} $ and  ${\mathcal{O}_{K_{+}}}^{*}$, respectively.  A character $ \chi $ on $ \Gamma(\mathfrak{f},\mathfrak{m}) $ is defined as follows: \[
\chi\left( \begin{pmatrix}
	a&b\\c&d
\end{pmatrix} \right) = \chi_0(a({\rm mod} ~\mathfrak{m}), d({\rm mod}~ \mathfrak{m}))\chi_1(ad-bc).
\]
Suppose $\textbf{k}:=(k_1,k_2,\cdots,k_n)\in(\mathbb{N}_{\geq 0})^n$, $\textbf{z}:=(z_1,z_2,\cdots,z_n)\in \mathbb{C}^n$, $ \gamma := (\gamma_1,\gamma_2,\ldots,\gamma_n)\in\Gamma $, for a subgroup $ \Gamma \subset GL_2^+(\mathbb{R})^n $ with $\gamma_{i} := \begin{pmatrix}
	a_{i}&{b_{i}}\\c_{i}&d_{i}
\end{pmatrix} =  \sigma_{i}  \begin{pmatrix}
	a & b\\c&d
\end{pmatrix} $ and $ a \in \mathbb{R}$. 
With all these, we can define a right weight $\textbf{k}$ action of $GL^{+}_{2}(K)$ on the space of complex-valued functions defined over $\mathbb{H}^{n}$ as follows:
$$
(f|_{\textbf{k}, \gamma }) (\textbf{z}) =  \displaystyle{\prod_{i=1}^{n}} \left[(det(\gamma_{i}))^{k_{i}/2} j(\gamma_{i}, z_{i})^{-k_{i}}\right] f(\gamma.\textbf{z}),
$$ where $ j(\gamma_{i}, z) = c_{i}z+d_{i}$,  for $ f:\mathbb{H}^n\to \mathbb{C} $, $\gamma \in GL_{2}^{+}(K)$ and $\textbf{z} \in \mathbb{H}^n$.

\smallskip

 A {\it classical Hilbert modular form} of weight $ \textbf{k} $ is a complex analytic function $ f:\mathbb{H}^n\to \mathbb{C} $ such that for each $ \gamma\in\Gamma $ and for each $ \textbf{z}\in\mathbb{H}^n $, it is invariant under the action of $PGL_{2}^{+}(K)= GL_{2}^{+}(K) / K^{*}$, i.e,   
\[
(f|_{\textbf{k}, \gamma }) (\textbf{z}) = \chi(\gamma) f(\textbf{z}),  
\]
and  $f$ is holomorphic at the cusps of $ \Gamma $. Such a function $ f $ also has a Fourier series expansion of the form 
\begin{equation}\label{FExp}
f(\textbf{z}) = \sum_\xi c(\xi) e^{2 \pi i Tr(\xi \textbf{z})},
\end{equation}
where $ Tr(\xi \textbf{z}) = \sum_{j=1}^n \sigma_j(\xi) z_j$. The sum  runs over zero and the totally positive elements of a lattice (depending on  $ \Gamma $). A Hilbert modular form $f$ is called as a {\it Hilbert cusp form} if it vanishes at each cusp of $\Gamma$, i.e., $ c(0)=0 $  in the Fourier expansion. Let $\mathcal{M}_{\textbf{k}}(\Gamma, \chi)$ ($\mathcal{S}_{\textbf{k}}(\Gamma, \chi)$) denote the space of classical Hilbert modular forms (respectively Hilbert cusp forms). These are finite-dimensional $\mathbb{C}$-vector spaces.

\subsubsection{Adelization of Hilbert Modular Forms}

Let  $ \mathfrak{M} $ denote  the ray class group associated to the modulus $ \mathcal{O}_{K}\mathcal{P}_\infty $ and $ h := |\mathfrak{M}|$. Let $ \{t_j\}_{j=1}^h $ be a collection of ideles with $ (t_j)_\infty = 1 $ for all $ 1\leq j\leq h $ such that $\{\mathfrak{t}_{j}\}_{j=1}^{h}$ forms a complete set of representatives for $ \mathfrak{M} $. We also let $ \mathfrak{t}_j $ be the fractional ideal of $ \mathcal{O}_F $, corresponding to the idele $ t_j $. For each $ j $ and for a fixed integral ideal $ \mathfrak{n} $, denote $ \Gamma_{j}{(\mathfrak{n})} $ as
the congruence subgroup $ \Gamma(\mathfrak{t}_j\mathfrak{d}, \mathfrak{n}) $ where $ \mathfrak{d}$ is the different ideal of $F$. For every $ 1\leq j\leq h $ let $ f_{j} $ be a classical Hilbert modular cusp form of fixed weight $ \textbf{k} $ for the congruent subgroup $ \Gamma_j(\mathfrak{n}) $ with the character $\chi_{0}: (\mathcal{O}_{K}/ \mathfrak{n})^{*} \rightarrow \overline{Q}^{*}$ with the following  Fourier expansion
\[ 
f_{j}(\textbf{z}) = \sum _{\substack{0\ll \mu \in \mathfrak{t}_j \\ }} a_{j}(\mu) e^{2 \pi iTr(\mu \textbf{z})}.
\]
 Shimura \cite{Shi} associates an $h$-tupple $(f_{1}, \ldots, f_{h})$ of classical Hilbert modular forms to an adelic automorphic form $F$ on $GL_{2}(\mathbb{A}_K)$, for $\mathbb{A}_K$ the adele ring of $K$. Take $\mathcal{S}_{\textbf{k}}(\mathfrak{n}, \chi_{0})$ as the space of cuspidal forms $\textbf{f} =(f_{1}, \ldots, f_{h}) $ where $f_{j} \in \mathcal{S}_{\textbf{k}}(\Gamma(\mathfrak{t}_j\mathfrak{d}, \mathfrak{n}), \chi_{0} )$. There exists a $j \in \{1, \ldots, h\}$ and a totally positive element $\mu \in \mathfrak{t}_j $ with $\mathfrak{n}= (\mu) {\mathfrak{t}_j}^{-1}$.  For any fractional ideal $\mathcal{\ell}$ of $K$, we define  
the Fourier coefficient $c(\mathcal{\ell}, \textbf{F})$ of $\textbf{F} \in \mathcal{S}_{\textbf{k}}(\mathfrak{n}, \chi_{0})$ as follows:
\begin{equation}
c(\mathcal{\ell}, \textbf{F})=
\begin{cases}
 a_{j}(\mu) N(\mathfrak{t}_{j})^{-k/2} \qquad  {\rm if}~  \mathcal{\ell} = (\mu) \mathfrak{t}_{j}^{-1}~ {\rm is} ~{\rm  integral}, \\
0 \qquad \qquad \qquad \qquad \qquad  {\rm otherwise}. \\
\end{cases}
\end{equation}
 Note that $c(\mathcal{\ell}, \textbf{F})$  depends neither on the choice of $\mathfrak{t}_j $'s nor on the choice of $\mu$.  Under some suitable conditions, there exists a finite order Hecke character (viewed as an idelic character) $ \Phi $ such that $ \textbf{F}(k g) = \Phi(k)\textbf{F}(g) $ for all $ k \in \mathbb{A}_K^{*} $ and $ g\in GL_2(\mathbb{A}_F) $. The space of such adelic Hilbert cusp forms of weight $ \textbf{k} $ and level $\mathfrak{n}$ and character $ \Phi $ produces a finite-dimensional $\mathbb{C}$-vector space and it is  denoted by $S_{\textbf{k}}(\mathfrak{n}, \Phi)$.  

\smallskip

The theory of  newforms associated to  Hilbert modular forms is mostly  developed by Shemanske and Walling (see \cite{Shi}, \cite{ST}). Newform theory allows  the decomposition of the space $S_{\textbf{k}}(\mathfrak{n},\Phi)$ into the space of old forms $S_{\textbf{k}}^{-}(\mathfrak{n},\Phi)$ and the space of newforms $S_{\textbf{k}}^{+}(\mathfrak{n},\Phi)$. A form $\textbf{F} \in S_{\textbf{k}}(\mathfrak{n},\Phi)$ is said to be a {\it newform} if $\textbf{F} \in S_{\textbf{k}}^{+}(\mathfrak{n},\Phi)$ and $\textbf{F}$ is a common eigenfunction for all Hecke operators $T_{\mathfrak{n}}(\wp)$, where $\wp$ is a prime ideal not dividing $\mathfrak{n}$. If $C(\mathcal{O}_{K},\mathbf{f}) =1$, then $\textbf{F}$ is said to be is a {\it normalized newform}.  Let $\textbf{F} \in S_{\textbf{k}}(\mathfrak{n},\Phi)$ be a normalised newform, then  eigenvalue of $\textbf{F}$ with respect to the Hecke operator $ T(\mathfrak{m}) $ is precisely the Fourier coefficient $C(\mathfrak{m}, \textbf{F})$.The  Fourier coefficient $C(\mathfrak{m}, \textbf{F})$ is a multiplicative function in $ \mathfrak{m} $, and it also satisfies Hecke type relation. Moreover, it satisfies a generalized conjecture of Ramanujan (see \cite{Bla}), i.e., 
 for every $ \epsilon>0 $, \[
C(\mathfrak{m}, \textbf{f}) \ll_\epsilon N(\mathfrak{m})^\epsilon. 
\]
 Associated to a normalized newform $\textbf{F}$, one can define the Dirichlet series 
\begin{equation*}\label{DSeries}
	L(s,\textbf{F}) = \sum_{\substack{0 \neq \mathfrak{m}\subset \mathcal{O}_K }} \frac{C(\mathfrak{m}, \textbf{F})}{N(\mathfrak{m})^s}. 
\end{equation*}
 It converges absolutely for $\Re(s)>1$ and  can be expressed as an Euler product also, valid  in $\Re(s)>1$ (similar to classical sense).

\subsubsection{Galois Representations} 
 Let $\textbf{F} \in S_{\textbf{k}}(\mathfrak{n},\Phi)$ be a newform with the Hecke eigenvalues $C(\mathfrak{m}, \textbf{F})$. Let $K_{\textbf{F}}$ denote the number field generated by the set of 
 eigenvalues $\{ C(\mathfrak{m}, \textbf{F})\}_{ \mathfrak{m}\subset \mathcal{O}_K }$ and $\mathcal{O}_{\textbf{K}}$ be its ring of integers. Let $\lambda$ be a prime of $\mathcal{O}_K$, and  $\mathcal{O}_{K,\lambda}$ be the completion of $\mathcal{O}_K$ at $\lambda$.
 
\smallskip

Wiles proved  the existence of Galois representations associated to certain Hilbert modular forms of parallel weight 	$k \ge 1$ attached to a totally real number field $K$. More precisely, he proved the following result. 

\begin{thm}\label{HMFRep}\cite[Theorem 1]{Wiles}
Let $\textbf{F} \in S_{\textbf{k}}(\mathfrak{n},\Phi)$ be a cuspidal newform of parallel weight $k \ge 1$. Let  $\textbf{F}$ be  ordinary at $\lambda$. Then there exist a continuous irreducible representation 
$$
\rho_{\textbf{F}, \lambda}: {\rm Gal}(\overline{K}/K) \rightarrow GL_{2}(\mathcal{O}_{K,\lambda})
$$
unramified outside $\mathfrak{n} N(\lambda)$, and such that, for all prime $\mathfrak{q} \nmid \mathfrak{n} N(\lambda)$, we have 
\begin{equation*}
Tr (\rho_{\textbf{F}, \lambda}) (\rm Frob_{\mathfrak{q}}) = C(\mathfrak{q}, \textbf{F}),\quad \mbox{and} \quad
det (\rho_{\textbf{F}, \lambda}) (\rm Frob_{\mathfrak{q}}) = \Phi(\mathfrak{q})N(\mathfrak{q})^{k-1}.
\end{equation*}
\end{thm}
The theorem follows from a result on the existence of $\Lambda$-adic representations
attached to $\Lambda$-adic modular forms (due to Hida for $K = \mathbb{Q}$ and to Wiles for geneal $K$), a lifting theorem of classical modular forms to $\Lambda$-adic modular forms (due to Hida for $k \ge 2$ and $K = \mathbb{Q}$, and to Wiles for $k \ge 1$ and general $K$), and the following theorem of Carayol \cite{Car86}. For any cuspidal newform of parallel weight $ k \geq 2$, Carayol proved the existence of continuous representation (as given in \thmref{HMFRep}) for any odd degree number fields and with some natural conditions for even degree number fields. 

\subsubsection{Certain applications:}

\noindent
 In this context, we again recall the aforementioned conjecture in a different way as  follows. 

 \smallskip
 
\begin{conjecture}(Eichler-Shimura)
Let $K$ be a totally real number field.  Let  $\mathfrak{f}$ be a Hilbert eigenform of parallel weight $2$ and level $\mathcal{N}$, with rational eigenvalues. Then, there is an elliptic curve $E_{\mathfrak{f}}$ over $K$ with conductor $\mathcal{N}$ having the same $L$-function as $\mathfrak{f}$.
\end{conjecture}

\smallskip

Towards this set-up,   Blasius \cite{Blasius} derived the following theorem, extending some techniques from the work of Hida \cite{Hida}.

\begin{thm}\label{BHida} (Blasius, Hida)
Let $K$ be a totally real field.  Let  $\mathfrak{f}$ be a Hilbert eigenform  of parallel weight $2$ and level $\mathcal{N}$, such that $\mathbb{Q}_{\mathfrak{f}} = \mathbb{Q}$. Suppose that \\
$(a)$ either degree of  $K / \mathbb{Q}$ is odd,\\
$(b)$ or there is a finite prime $\mathfrak{q}$ such that  $\pi_{\mathfrak{q}}$ belongs to the discrete series, where $\pi$ is the cuspidal  automorphic representation of $GL_{2}(\mathbb{A}_{K})$ attached to $\mathfrak{f}$.\\
Then, there exists an elliptic curve $E_{\mathfrak{f}} $ over $K$ of  conductor $\mathcal{N}$ having the same $L$-function as $\mathfrak{f}$.
\end{thm}
 If $v_{\mathfrak{q}}(\mathcal{N})=1$
then $\pi_{\mathfrak{q}}$ is special and $(b)$ is satisfied. This is the only case of $(b)$ that is required for the  modularity of $E$ over real fields . Proofs of this case of \thmref{BHida} is given by Darmon \cite{Darmon} and Zhang \cite{Zhang}, independently. We end this section by mentioning the following corollary which deals about the existence of an elliptic curve with same L-function as the considered Hilbert eigenform.

\begin{cor}\label{Cor1}
Let $E$ is an elliptic curve over a totally real field $K$ and $p$ be an odd rational prime. Suppose  $\overline{\rho}_{E,p}$ is irreducible, and  $\overline{\rho}_{E,p}  \sim \overline{\rho}_{\mathfrak{f},\varpi}$ for some Hilbert eigenform (newform) $\mathfrak{f}$ of parallel weight $2$ with rational eigenvalues.  Let $\mathfrak{q}\nmid p$ be a prime in $K$ such that the following conditions hold:
\begin{enumerate}
\item $E$ has potentially multiplicative reduction at $\mathfrak{q}$,
\item $p\mid \# \overline{\rho}_{E,p}(I_{\mathfrak{q}})$.
\item $p \nmid ({\rm Norm_{K/\mathbb{Q}}}\pm 1)$.  
\end{enumerate}
Then, there exists an elliptic curve  $E_{\mathfrak{f}}/K$ of  conductor $\mathcal{N}$ with the same
$L$-function as $\mathfrak{f}$.
\end{cor}

 
 \section{Elliptic curve over  certain real  fields}

In this direction, beyond the field of rational numbers, the first breakthrough came for real quadratic fields. It has been  generalized later for cubic fields. Therefore, in this exposition, we  try to  give an overview of the way from $\mathbb{Q}$ to real quadratic fields. We also emphasize understanding the difference (hardness)  between these two fields. We end this section by briefing the route towards cubic fields. 

\subsection{Elliptic curve over real quadratic fields}
Let  $E$ be an elliptic curve over a totally real number field $K$. Due to the previous discussion, it is clear to us that for understanding elliptic curves over totally real fields, we need to understand the Galois representations, namely $\overline{\rho}_{E,p} (G_{K(\zeta_p)})$. In that regard,   the following theorem of Breuil and Diamond \cite{BD2014} tells us about the exact properties of $\overline{\rho}_{E,p} (G_{K(\zeta_p)})$ which we need to consider. More precisely, it connects the modularity of the elliptic curve and the crucial conditions on the image of $\overline{\rho}_{E,p} (G_{K(\zeta_p)})$.
\begin{thm}\label{2 conditionquadratic}
Let  $E$ be an elliptic curve over a totally real number field $K$ and let $ p \neq 2$ be a rational prime.  Suppose 
\begin{itemize}
\item $\overline{\rho}_{E,p}$ is modular
\item $\overline{\rho}_{E,p} (G_{K(\zeta_p)})$ is absolutely irreducible.
\end{itemize} Then $E$ is modular.
\end{thm}

Along with this theorem, we also have the aforementioned Tunnell's theorem ($ \overline{\rho}_{E, 3}$  is irreducible $\Rightarrow$ $ \overline{\rho}_{E, 3}$ is modular) and $3-5$ modularity switching argument. Due to these, one can deduce the following important theorem. 

    \begin{thm} \label{3-5 switching} \cite{FHS2015}
 Let $E$ be an elliptic curve over a totally real field $K$ and $ p = 3$ or $5$. If $\overline{\rho}_{E,p}(G_{K(\zeta_p)})$ is absolutely irreducible, then $E$ is modular.
\end{thm}

Now get back to the steps of the proof of the above theorem over real quadratic fields. Due to the theorem of Langlands and Tunnell \cite{Langlands, Tunnell}, if $E$ is an elliptic curve over a totally real field $K$ and $\overline{\rho}_{E,3}(G_{K(\zeta_3)})$ is absolutely irreducible, then  $\overline{\rho}_{E,3}$ is moduler. Therefore, using Theorem \ref{2 conditionquadratic}, we have $E$  is modular. In the remaining part of this modularity case, the concept of $ 3$ - $ 5$ modularity switching and $ 3$ - $7 $ modularity switching have been used.  To be more precise, one has to show that there exists an elliptic curve $ E ^ \prime$ over the same real quadratic field $K$ such that $$ \overline{\rho}_{E, 5} \sim \overline{\rho}_{E^\prime, 5} $$ along with the condition that $ \overline{\rho}_{E^\prime, 3}(G_K)$ contains $ SL_2(\mathbb{F}_3)$. Using above discussion, we get that $ E^\prime$ is modular. Now, by the equivalence notation above, we can derive that $ \overline{\rho}_{E,5}$ is also modular. Therefore, all conditions in Theorem \ref{2 conditionquadratic} with $p = 5$ have been satisfied. So, we conclude that $ E$ is modular.  For applying $3$-$5$ modularity switching, it is  enough to find out an Elliptic curve $ E^\prime/K$ with  equivalence Galois representation. With all these discussion, Theorem \ref{3-5 switching} follows. 

\smallskip

Due to Theorem \ref{3-5 switching},  it is enough to parametrize elliptic curves where $\overline{\rho}(G_{K(\zeta_3)})$ and $\overline{\rho}(G_{K(\zeta_5)})$ are simultaneously absolutely reducible. Some of these modular curves  have infinitely many real  quadratic points and there seems to be no easy way to prove that these points are modular. For that purpose, $3-7$ modularity switching argument plays a crucial role to overcome this obstacle. 
\begin{thm} \label{3-7 switching} [$3-7$ modularity switching argument]
Let $E$ be an elliptic curve over a totally real field $K$. If $ \overline{\rho}_{E,7}(G_{K(\zeta_7)})$ is absolutely irreducible, then $E$ is modular.
\end{thm}

When the elliptic curve $E$ is defined over the real quadratic field $K$, then in \cite{FHS2015}, following path was considered. They have proved (in Lemma 3.2 \cite{FHS2015}) that there is a totally real quartic extension $ L$ of $K$ and an elliptic curve $ E^\prime$ over $ L $ such that 
\begin{itemize}
\item $ L \cap K(E[7]) = K$ 
\item $\overline{\rho}_{E,7} \mid G_7 \sim \overline{\rho}_{E^\prime, 7}$ and $ \overline{\rho}_{E^\prime, 3} (G_L)$ contains $ SL_2(\mathbb{F}_3)$.
\end{itemize} Since $ \overline{\rho}_{E, 7}(G_{K(\zeta_7)})$ is  absolutely irreducible, by the above discussion we have observed that $ E ^ \prime$ is modular.  Therefore $ \overline{\rho}_{E^\prime, 7}$ is modular which implies $\overline{\rho}_{E,7} \mid G_7$ is also modular. Assume that $\overline{\rho}_{E,7} (G_{L(\zeta_7)})$ is absolutely irreducible.  Again by Theorem \ref{2 conditionquadratic}, one can deduce that $ E/L$ is modular. By repeated application of cyclic base change  results of Langlands \cite{Langlands, Tunnell}, it follows that $ E / K $ is modular.  It remains to show that $\overline{\rho}_{E,7} (G_{L(\zeta_7)})$ is absolutely irreducible.  To obtain that we can use the facts that $\overline{\rho}_{E,7} (G_{K(\zeta_7)})$ is absolutely irreducible and  the above mentioned fact $ L \cap K(E[7]) = K$.  This concludes the sketch of  the proof of Theorem \ref{3-7 switching}.

 \bigskip

From the above observations, it is clear that if an elliptic curve $E$ over a totally real field $K$ is modular except if the images $\overline{\rho}(G_{K(\zeta_p)})$ are simultaneously absolutely reducible for $p = 3, 5, 7$. In other words,  one needs to understand the second condition of Theorem \ref{2 conditionquadratic} more closely. In \cite[Proposition 12]{FHS2015}, it has been shown that  whenever any elliptic curve violates the aforementioned second condition, then   $E$ gives rise to non-cuspidal $K$-points on certain modular curves. Therefore, it is clear that we need to examine some modular curves carefully and tackle this.  For that purpose, we recall the modular curves here briefly.  

\smallskip

\noindent {\bf Modular curves:} 
Let $N$ be a positive integer and $ G$ be a subgroup of $ GL_2(\mathbb{Z}/ N \mathbb{Z})$. We know $\det(G) \subset (\mathbb{Z}/ N \mathbb{Z})^{\times}$. Let $ \zeta_N$ be a primitive $N$-th root of unity. By the action of $ \det(G) $ on $ \zeta_N$, we get a fixed subring of $ \mathbb{Z} [ 1/ N, \zeta_N]$, say $ R_G$. We can associate  $G$ to a {\it modular curve } $ X_G$ (for more details \cite{KMBook}). 

\smallskip

Let $ \mathbb{H}^{*} = \mathbb{H} \cup \mathbb{Q} \cup \{ \infty\}$ be the extended upper half plane.  The full modular group $ SL_2(\mathbb{Z})$ acts on $ \mathbb{H}^{*}$ by the usual fractional linear transformations.  The quotient $ SL_2(\mathbb{Z})  \setminus \mathbb{H}^*$  becomes a compact Riemann surface of genus $0$ and therefore  an analytic structure can be given on it over $ \mathbb{C}$  which we denote  by $ X(1)$. The  set $ \mathbb{Q} \cup \{ \infty \}$  forms a single orbit under the action of $ SL_2(\mathbb{Z})$ such that the orbit corresponds to a point on the curve $ X(1)$ which we call  cusp $\infty$.  Let $ H$ be a subgroup of $ GL_2(\mathbb{Z}/ p \mathbb{Z})$ satisfying $ {\rm det}(H) =(\mathbb{Z}/ p \mathbb{Z})^{\times}$.   We associate a congruence subgroup $(\Gamma_H)$ to $ H$ as  the preimage of $ H \cap SL_2(\mathbb{Z}/ p \mathbb{Z})$ under the map $ SL_2(\mathbb{Z}) \rightarrow SL_2(\mathbb{Z}/ p \mathbb{Z})$.   The modular curve $ X_H / \mathbb{C}$ is the algebraic  curve whose analytic version is   the compact Riemann surface $ \Gamma_H \backslash \mathbb{H}^* $. The inclusion $ \Gamma_H \subset SL_2 (\mathbb{Z})$ includes a morphism $ j : X_H \rightarrow X(1)$  (which is also defined over Spec$(\mathbb{Z}[1/p])$ as discussed above. The cusps of $ X_H$ are the pre-images of $ \infty$  in $X_{H}$ and therefore the orbits of $ \mathbb{Q} \cup \{ \infty \}$ under the action of $ \Gamma_H$. 

\smallskip

 Now we briefly mention the modular interpretation of rational points on $ X_H$. Furthermore, we assume that $ -I \in H$. Let $K $ be a field of characteristic $ \neq p$. 
 \begin{enumerate}
     \item[$\bullet$] Let $ E $ be an elliptic curve over $ K$ such that $ \overline{ \rho}_{ E, p}(G_K)$ is conjugate to a subgroup of $ H$. There is at least one noncuspidal point $ P \in X_H(K)$ such that $ j(P) -= j(E)$.
     \item[$\bullet$] There is a partial converse available. It says that if $ P \in X_H(K)$ is a noncuspidal point and $ j(P) \neq 0, 1728$, then there is an elliptic curve $ E$ over $ K$ such that $ \overline{\rho}_{E,p}(G_K)$ is conjugate to a subgroup of $ H$ and $ j(H) = j (P)$.
 \end{enumerate}
 If we consider $ H = B_0(p)$ (Borel), $C_s^+(p)$ (split Cartan), $C_{ns}^+(p)$ (non-split Cartan) subgroups of  $GL_2(\mathbb{F}_p)$, then corresponding $ X_H$ is the modular curve usually denoted by $ X_0(P)$, $ X_{split}(P)$ and $ X_{nonsplit}(P)$ respectively.

\smallskip

 It have been shown (in  \cite{FHS2015}) that  if  $\overline{\rho}_{E,p}(G_{K(\zeta_p)})$ is reducible, then the image of  $ \overline{\rho}_{E,p}$ lies either in a Borel subgroup of $ GL_2( \mathbb{F}_p)$ or in certain subgroups contained in the normalizers of Cartan subgroups.

\smallskip

 Now we pass to the more crucial case, where one needs to consider the cases where $ \overline{\rho}_{E,p}( G_{K(\zeta_P)})$  for $ p =3, 5$ and $7$ are absolutely reducible. As we discussed earlier, in those cases, the elliptic curve  $ E/K$ produces a non-cuspidal $K$-point  on a modular curve corresponding to one of these subgroups.  In \cite{FHS2015}, it has been observed that such an elliptic curve gives rise to a $K$-point on one of $27$ modular curves. In \cite{FHS2015}, Freitas et al. tackeled each of these cases using the theory of modular forms. At first, it has been proved that any elliptic curve over $ \mathbb{Q}(\sqrt{5})$ are modular. For further, one can start with that $ K \neq \mathbb{Q}(\sqrt{5})$. Then  it has been proved that {\it all non-cuspidal real quadratic points on the  following modular curves correspond to  modular elliptic curves. 
 \begin{equation} \label{important list of modular curves}
 X(b5, b7),  X(b3, s5), X(s3, s5), X(b3, b5, d7), X(s3, b5, d7), X(b3, b5, e7), X(s3, b5, e7),
 \end{equation} where $e7$, $d7$ indicate mod-$7$ level structures that are finer than normalizer of split-Cartan and non-split Cartan subgroups.  b, s and ns  stand for {\it Borel}, {\it normalizer of split Cartan} and {\it non-split Cartan}} respectively. To  explain more clearly, non-cuspidal $K$-points on $ X(s3, b5, d7)$ provide elliptic curves $ E$ over $K$ such that the image of $ \overline{\rho}_{E,3}$ is contained in a conjugate of the normalizer of a split Cartan subgroup of $ GL_2(\mathbb{F}_3)$, the image of  $ \overline{\rho}_{E,5}$ is contained in a conjugate of the Borel subgroup of $ GL_2(\mathbb{F}_5)$ and the image of $ \overline{\rho}_{E,7}$ is contained in a conjugate of some group of $ GL_2(\mathbb{F}_7)$. For more details, we record the following result which details about the image of $ \overline{\rho}_{E,p}$, when $\overline{\rho}_{E, p} (G_{K(\zeta_p)})$ is absolutely reducible (for $p = 3$ or $5$ or $7$). 
\begin{thm} \cite[Proposition 2.1]{FHS2015}
 Let $E$ be an elliptic curve over a totally real field $K$. Let $ p \geq 3$ be a rational  prime. Then 
 \begin{itemize}
 \item $\overline{\rho}_{E, p} (G_{K(\zeta_p)}) = \overline{\rho}_{E,p}(G_K) \cap  SL_2(\mathbb{F}_p)$.
 \item If $ \overline{\rho}(G_{K(\zeta_p)})$  is absolutely reducible, then $\overline{\rho}_{E, p} (G_K)$ is contained either in a Borel subgroup or in the normalizer of a Cartan subgroup. In this case $ E$  generated to a non-cuspidal $ K $ -point on $ X_0(p)$ , $ X_{\text{split}}(p)$ or $ X_{\text{nonsplit}}(p)$; where $X_0(p)$, $X_{\text{split}}(p)$ and $X_{\text{nonsplit}}(p)$ signify modular curves corresponding to  Borel subgroup :$B(p)$, split-Cartan: $ C_s(p)$ and non-split cartan: $ C_{ns}(p)$, respectively.  
 \end{itemize}
\end{thm}

 Now get back to equation  \eqref{important list of modular curves}. In the literature, there are explicit methods which in certain special cases are capable of determining all quadratic points on a given curve. Freitas et al. (in \cite{FHS2015}) cultivated the known literature in this area further to obtain the quadratic points on the three modular curves, $X(b5, b7),  X(b3, s5), X(s3, s5)$,  and deduce modularity of the non-cuspidal real quadratic points.  The other four cases in the above list are complicated because of their large genera.  In \cite{FHS2015}, representations of these curves have been given as normalizers of fibre products of  suitable curves. By some nice trick, they could able to show that the non-cuspidal  real quadratic points on these curves have $ j$-invariants belonging to $ \mathbb{Q}$ and thus our well-known modularity over rational field can be applied.  This is the way, we get an overview of the following result. 
 \begin{thm}\cite{FHS2015}
     Let $E$ be an elliptic cuve over a real quadratic field $K$. Then we have $E/K$ is modular. 
 \end{thm}
 Next, we  emphasize the main steps of the proof of modularity of elliptic curves over real quadratic fields by a diagram. 
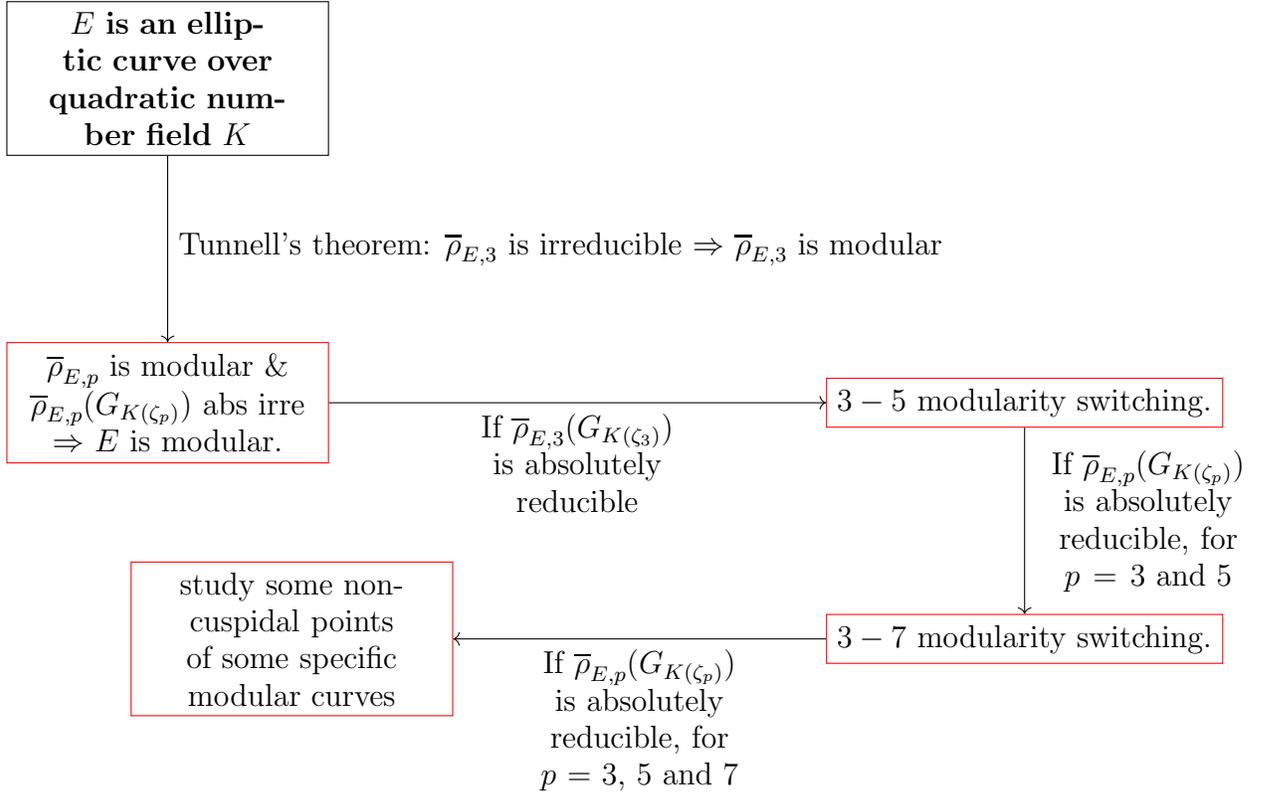
\begin{center}
        \begin{tikzpicture}[node distance=3.0cm]
\node(a)[process] {{\textbf{$E$ is an elliptic curve over quadratic  number field $K$}}};
\node(b)[process, draw=red, below= 6 em of a] {$\overline{\rho}_{E, p}$ is modular 
 $\&$ $\overline{\rho}_{E, p}(G_{K(\zeta_p)})$ abs irre $\Rightarrow$  $E$ is modular.};
\node(c)[rectangle, draw=red, right=16em of b] {$3-5$ modularity switching.};
\node(d)[rectangle, draw=red, below=6em of c] {$3-7$ modularity switching.};
\node(e)[process, draw=red, left=12em of d] { study some non-cuspidal points of some specific modular curves};
\draw[->] (a) -- node[distance=14.0cm, right] {Tunnell's theorem: $ \overline{\rho}_{E, 3}$  is irreducible $\Rightarrow$ $ \overline{\rho}_{E, 3}$ is modular}(b);
\draw[->] (b) -- node[process1, distance=3.0cm, below] 
 {If $\overline{\rho}_{E, 3}(G_{K(\zeta_3)})$ is absolutely reducible}(c);
\draw[->] (c) -- node[process1, distance=3.0cm, right] {If $\overline{\rho}_{E, p}(G_{K(\zeta_p)})$ is absolutely reducible, for $p=3$ and $5$}(d);
\draw[->] (d) -- node[process1, distance=3.0cm, below] {If $\overline{\rho}_{E, p}(G_{K(\zeta_p)})$ is absolutely reducible, for $p=3$, $5$ and $7$}(e);
\end{tikzpicture}
\end{center}
 \bigskip

 \subsection{Elliptic curves over real cubic fields}
 Now we can move to the cubic field case.  By  a result of Rubin \cite[Proposition 6]{rubin1997}, we know that if $ \overline{\rho}_{E,3}(G_{K(\zeta_3)})$  is absolutely  reducible, then it is conjugate to a subgroup of  $B(3)$ or  $ C_s^+(3)$. For notation purposes, we denote $B(p)$ as the Borel subgroup of $ GL_2(\mathbb{F}_p)$ and $ C_s^+(p)$ denoted the normalizer of the split Cartan subgroup $ C_s(p)$ as 
 we did earlier. Combining with Theorem \ref{3-5 switching}, we know that  if $ \overline{\rho}_{E,3}(Gal (\overline{K}/K))$ is not conjugate to a subgroup of  $ B(3)$ or $ C_s^+(3)$, then $ E$ will be modular. 

 \smallskip
 
 Now it is time to consider $p=5$ and $ p=7$ as we have  noticed in the above quadratic case.  After the modularity theorem of totally real quadratic field  appeared in the literature,  Throne \cite{thorne2016} considered $p=5$ case more deeply in the following theorem which will reduce some hurdles as mentioned above. 
 
 \begin{thm}\cite{thorne2016}
 Let $K$ be a totally real field and $E$ be an elliptic curve over $K$. Suppose $5$ is not a square in $K$. Moreover, if $\overline{\rho}_{E,5}$ is irreducible, then we get  $E$ is modular. 
 \end{thm} Also Kalyanswamy \cite{kal2018} simplified the case $p=7$ in the following sense. 
 
 \begin{thm}\cite{kal2018}
 Let $ E$ be an elliptic curve over a totally real field $K$. If the following conditions hold.
 \begin{enumerate}
 \item $ K \cap \mathbb{Q}(\zeta_7) = \mathbb{Q}$
 \item $ \overline{\rho}_{E,7}$ is irreducible. 
 \item $ \overline{\rho}_{E, 7}(Gal (\overline{K}/K))$ is not conjugate to a subgroup of $ C_{ns}^+(7)$.
 \end{enumerate}
then $ E$ is modular. 
 \end{thm} 
 
 Note that $ K = \mathbb{Q}(\zeta_7)^+$ is the only  real cubic field for which the above theorem is not applicable.  Now we  consider that $ K $ is a totally real cubic field. For the time being, we also consider $ K \neq 
 \mathbb{Q}(\zeta_7)^+$ .  It is clear from the above discussion that if $ E$ is non-modular over $ K$, then  it would rise to a non-cuspidal $K$-point $P$ on  either $ X_0(3)$ or on $ X_{\text{split}}(3)$. For $ p=5$ or $7$, simultaneously,  it would give a non-cuspidal $K$-point $ Q$ on $X_0(5)$ and similarly a non-cuspidal  $ K$ point $ R $ on either $ X_0(7)$ or $ X_{\text{nonsplit}}(7)$.  Hence we obtain a $ K$ -point on one of the fiber products $$  X_u(3) \times _{X(1)} X_0(5) \times_{X(1)} X_v(7) \qquad u \in \{0, \text{split}\}, v \in \{ 0, \text{nonsplit}\} $$  Denoting the normalization of the above as $ X(u3, 05, v7)$, we observe that when $E$ is assumed to be non-modular, then it is also non-CM. Therefore, $ j(E)$ can not be $0$, $1728$.  The maps $ X_H \rightarrow X(1)$ are ramified only above $0$, $1728$ and $ \infty$ and hence the $ K$-point we obtain from $ E$ on the above product  is a smooth point and hence gives rise to a $ K$ -point on the normalization $ X(u3, 05, v7)$. \\
 
For concluding  modularity of elliptic curves over totally real cubic fields, it is enough to prove that $ K$ -points on the four possible curves $ X(u3, 05, v7)$ as cuspidal. Therefore, it is enough to prove the following two cases. 
\begin{thm}
 \begin{enumerate}
\item Let $K \neq \mathbb{Q}(\zeta_7)^+$ be a totally  real cubic field. Then the curve  $ X(05, 07)(K)$ consists only of cusps. 
\item Let $K$ be a cubic field. Then $ X(05, \text{nonsplit }7)(K)$ consists only of cusps.
\end{enumerate}
\end{thm}
In the paper \cite{DNS2020}, Derickx et al. proved the aforementioned two claims to obtain the modularity of elliptic curves over totally real cubic fields. Moreover,  they tackled the case  when $ E$ is  an elliptic curve defined over $ K = \mathbb{Q}(\zeta_7)^+$. The technique was to understand the  modular curves and the corresponding jacobians. In this way, we can briefly understand the steps of the proof of the following theorem. 
\begin{thm}
    Let $K$ be a totally real cubic number field and $ E$ be an elliptic curve over $K$. Then $E$ is modular. 
\end{thm}

\bigskip

 \section{Applications of Modularity theorems}
 The {\it modularity  of elliptic curves} is connected with the prrof of  the famous and  long-standing conjecture ``Fermat's Last Theorem" which says that the equation  $$ x^{n} +y^{n} = z^{n}  \quad {\rm with~ } x,y,z \in \mathbb{Z}$$ 
has no non-trivial solution (means $xyz \not = 0$) over $\mathbb{Z}$ for $n \ge 3$.  The modularity theorem  of elliptic curves over $\mathbb{Q}$ implies  ``Fermat Last Theorem"  over $\mathbb{Q}$. This was proved by Wiles \cite{Wiles}  and Taylor-Wiles \cite{Taylor-Wiles}. Extending this concept,  Debarre-Klassen \cite{Deb-Kl}   made the following conjecture.

\smallskip
\noindent
\textbf{Conjecture:} Let $K$ be a number field of degree $d$ over $\mathbb{Q}$. Then the equation $x^{n} +y^{n} = z^{n}$ has only trivial solution over $K$ when $n \ge d+2$.

\noindent
 Here, the trivial solution means to a point $(a,b,c)$ on the curve $x^{n} +y^{n} = z^{n}$ such that $a+b=c$. This conjecture deals not only with rational points but also with the solutions of the form $\omega^n+\overline{\omega}^n=1$ where $\omega$ is the primitive $6^{\rm th}$ root of unity, belonging to any fields containing $\mathbb{Q}(\sqrt{-3})$ with  $n \equiv 1$ or $5 (\rm mod ~6)$.

As we mentioned earlier the first result towards this direction was done  over $\mathbb{Q}(\sqrt{2})$  in \cite{JM2004} and \cite{JM2008}. More Precisely,  the following result was obtained first, later it was proved for any semistable elliptic curves. 

\begin{thm}\label{ModularyQ(2)} [assuming an unpublished result of Fujiwara \cite{Fujiwara}]
The equation $ x^{n} +y^{n} = z^{n}  ~ {\rm with }~ x,y,z \in \mathbb{Z}[\sqrt{2}]$ has no solution with $xyz \not =0$ when $n \ge 4$.
\end{thm}

The above result was obtained by extending the work of Ribet \cite{Ribet} and Wiles \cite{Wiles}. It was also observed that there are no other real quadratic fields $\mathbb{Q}(\sqrt{d})$  for which a similar argument of Ribet and Wiles's method can be applied. In \cite{FS2014}, Freitas-Siksek obtained that there are no non-trivial solutions to Fermat's Last Theorem over $\mathbb{Q}(\sqrt{d})$ where $d \not = 5, 17$ is a square-free integer $3 \le d \le 23$ with $n \ge 4$.  In \cite{PMJ22}, Michaud-Jacobs extended  for most square-free  $d$ in the range $26 \le d \le 97$ with  $n \ge 4$. 

\smallskip

In general, for tackling Fermat's Last Theorem over real quadratic fields of a larger discriminant,  the obstacle has encountered that  the  demonstration of the irreducibility of certain Galois representations and elimination of the number of Hilbert newforms. To overcome such obstacles, Freitas and Siksek \cite{FS2015} come up with the idea of studying  Fermat's Last theorem asymptotically, namely {\it Asymptotic Fermat’s Last Theorem} for certain families of number fields. More precisely, they state that {\it there exist  an  absolute constant  $B_{K}$ depending on the totally real field $K$ such that  the Fermat's equation $ x^{p} +y^{p} = z^{p}, ~ x,y,z \in K$ with $p > B_{K}$ has only trivial solution.} If the absolute constant is effectively computable, then we shall refer to this as the {\it effective asymptotic Fermat’s Last Theorem over $K$}. To  address this version of Fermat's last theorem,  the following theorem  was proved in \cite{FS2015}.
\begin{thm} (Freitas and Siksek \cite{FS2015})
Let $d \ge 2$ be a square-free integer, satisfying one of the following conditions:
\begin{enumerate}
\item $d \equiv 3 ({\rm mod~8})$.
\item $d \equiv 6$ or $10 ({\rm mod~8})$.
\item$ d \equiv 2 ({\rm mod~16})$ and d has some prime divisor $q \equiv 5 {\rm ~or~} 7 ({\rm mod~8})$.
\item $d \equiv 14 ({\rm mod~16})$ and d has some prime divisor $q \equiv 3 {\rm ~or~} 5 ({\rm mod~8})$.
\end{enumerate}
Then, the  effective asymptotic Fermat’s Last Theorem holds over $K =\mathbb{Q}(\sqrt{d})$.
\end{thm}

With the help of modularity, level lowering (see for detail; \cite{Fujiwara, Jarvis, Raj}) and image of inertia comparisons ( see; \cite{bernettskinner}), they give an algorithmically testable criterion which (if satisfied by $K$) implies the asymptotic Fermat’s Last Theorem holds over $K$. Using techniques from analytic number theory, they show that the
criterion is satisfied by $K = \mathbb{Q}(\sqrt{d})$ for a subset of $d \ge 2$ having density $5/6$
among the squarefree positive integers. One can improve this density to 1 if 
 a standard “Eichler–Shimura” conjecture is assumed to be true.

In recent work \cite{FKS2020},  Freitas-Kraus-Siksek  study the solutions to certain S-unit equations to show that Asymptotic Fermat’s Last Theorem holds for several families of number fields, including some real biquadratic fields.

\begin{thm}  (Freitas-Kraus-Siksek, \cite{FKS2020})

The effective Fermat’s Last Theorem holds over $K =\mathbb{Q}(\sqrt{d}, \sqrt{\ell})$ for prime $\ell \equiv 3  ({\rm mod~8})$.
\end{thm}

 It is a natural problem then to study Fermat's Last Theorem over real biquadratic fields such as real biquadratic fields $\mathbb{Q}(\sqrt{2}, \sqrt{3})$. Recently,  Khwaja and Jarvis \cite{MKFJ22}  prove the following result and discuss some obstacles that arise over more general real biquadratic fields.
\begin{thm} (Khwaja and Jarvis \cite{MKFJ22})
Let $K= \mathbb{Q}(\sqrt{2}, \sqrt{3})$ and $n \ge 4$ . Then, the Fermat equation $x^{n} +y^{n} = z^{n}$  has no non-trival solution over $\mathbb{Q}(\sqrt{2}, \sqrt{3})$.
\end{thm}
\noindent
 There are some obstacles that arise when extending the methods to more general real biquadratic fields.
Moreover, they also prove the following result.

 \begin{thm} (Khwaja and Jarvis \cite{MKFJ22})
Let $K= \mathbb{Q}(\sqrt{2}, \sqrt{11})$, and $p \ge 29$ and $p \not= 197$.  The Fermat equation $x^{p} +y^{p} = z^{p}$  has no non-trival solution over $\mathbb{Q}(\sqrt{2}, \sqrt{11})$ for $n \ge 4$.
\end{thm}
The exclusion of prime $p = 197$ arises during the application of a result of Kraus (in \cite{Kraus}) while proving the irreducibility of $\overline{\rho}_{E,p }$. 
 
\smallskip

 Generalizing the aforementioned concept, it is natural to  examine the solutions of {\it generalised Fermat's equations} which is  denoted as  $$ A x^{p} +B y^{q} + Cz^{r} =0,$$ where $p, q, r$ are natural numbers.  Ribet \cite[Theorem 3]{Ribet97} proved that there are no non-trivial integral solution to   $  x^{p} +2^{r} y^{p} + z^{p} =0$ with $1 \le r\le p$ over $\mathbb{Q}$. Darmon and Merel \cite{Darmon} established that $  x^{n} + y^{n}  = 2 z^{n}$ has no non-trivial primitive integral  solution over $\mathbb{Q}$ for $n \ge 3$. Many more results in this direction has been established. Following the idea of Freitas-Siksek \cite{FS2015}, Deconinck \cite[Theorem 1]{Dec16} studied asymptotic version of  generalised Fermat's Last theorem $ A x^{p} +B y^{p} = Cz^{p}$ with $2 \nmid ABC$, and established the asymptotic result. In his work, Kumar-Sahoo \cite{NKSS2022}, studied the following  Diophantine equations  $  x^{p} + y^{p} = 2^{r} z^{p}$ of exponent $p$ and $r \in \mathbb{N}$, and  $  x^{p} + y^{p} =  z^{2}$ over totally real fields. Precisely, they established the following result.
 
\begin{thm} (Kumar-Sahoo \cite[Themem3.2]{NKSS2022} )
Let $K$ be a totally real field and $S_{K} = \{ \mathcal{P} : \mathcal{P} ~is ~ a ~ prime ~ ideal ~of~ \mathcal{O}_{K}~ with ~ \mathcal{P} \mid 2\}$. Suppose for every $(\lambda, \mu)$ to the $S_{K}$-unit solution of 
$$
\lambda + \mu =1, ~~ \lambda, \mu \in \mathcal{O}^{*}_{S_K},
$$ 
there exist some $\mathcal{P} \in S_{K}$ that satisfies 
$$
{\rm max} (|v_{\mathcal{P}}(\lambda)|, |v_{\mathcal{P}}(\mu)|) \le 4v_{\mathcal{P}}(2). 
$$ 
Then, there exists a constant $B_{K,r}$ such that for prime $p > B_{K,r}$, the  Diophantine equations  $  x^{p} + y^{p} = 2^{r} z^{p}$ of exponent $p$ has no non-trivial solution in the ring of integers of $K$.
\end{thm} 
  In \cite{IKO2020}, Isik et al. proved that the Diophantine equation  $  x^{p} + y^{p}  =  z^{2}$ does not have any non-trivial solutions of certain type over $K$, if there exists a prime $\mathcal{P}$ (dividing $2$) in $K$ with inertial degree $f(\mathcal{P},2)=1$. This result has been established by Kumar-Sahoo without any assumption on the prime $\mathcal{P}$ in $K$. They also provide several local criterion of $K$ for which the Diophantine equations $  x^{p} + y^{p} = 2^{r} z^{p}$ of exponent $p$ and $r \in \mathbb{N}$ and  $  x^{p} + y^{p} = 2^{r} z^{p}$ with $r=2,3$ has no non-trivial solution of certain   types over $K$.

      
  \footnotesize{
  	 
}


\begin{thebibliography}{10}  
  	

  	

 

   		
 
 
 
 
  
\bibitem{bernettskinner}
M. A.  Bennett, Michael and C. M. Skinner, 
{\it Ternary Diophantine equations via Galois representations and modular forms}, Canad. J. Math. 56 (2004), no. 1, 23–54. 23–54.

\bibitem {Bla} D. Blasius;
 	{\it Hilbert modular forms and Ramanujan conjecture},
 	Noncommutative geometry and Number Theory, Wiesbaden, in: Aspects Math. vol. E37, Vieweg, 2006, pp. 35-56.
  
  \bibitem{Blasius}
  D. Blasius, {\it Elliptic curves, Hilbert modular forms, and the Hodge conjecture, Contributions
to automorphic forms, geometry, and number theory}, 83--103, Johns Hopkins Univ. Press,
2004.
 
 
 
 
 \bibitem{Box2021}
 Josha Box, {\it Elliptic curve over totally real quartic field not containing $\sqrt{5}$ are modular}, arXiv Pre-print (2021) 
 
 
  \bibitem{BCDT} 
 C. Breuil, B. Conrad, F. Diamond and R. Taylor, {\it On the modularity of elliptic curves over Q: wild 3-adic exercises}, Journal of the American Mathematical Society {\bf 14} (2001), 843--939.
 
 
 \bibitem{BD2014} 
 C. Breuil and F. Diamond, {\it Formes modulaires de Hilbert modulo p et valeurs d’extensions galoisiennes}, Annales Scientifiques de l’Ecole Normale Sup´erieure ´ {\bf 47} (2014), no. 5, 905--974.
 \bibitem{Bruin2000}
N. Bruin, On powers as sums of two cubes. Algorithmic number theory (Leiden, 2000), volume 1838
of Lecture Notes in Comput. Sci., pp. 169–184. Springer, Berlin (2000)

 \bibitem{Car86}
H. Carayol, { \it Sur les repr\'{e}sentations $\ell$-adiques associ\'{e}es aux formes modulaires de Hilbert}, Ann. Sci. Ec. Norm. Sup.  {\bf 19}
(1986), 409 --468.
 
\bibitem{CS2009}
I. Chen, S. Siksek, Perfect powers expressible as sums of two cubes. J. Algebra 322(3), 638–656
(2009)
 
\bibitem{David}
A. David,  {\it Caract\'{e}re d’isog\'{e}nie et crit\'{e}res d’irr\'{e}ductibilit\'{e}}, (2011). https://arxiv.org/abs/1103.3892

\bibitem{Darmon}
H. Darmon, {\it Rational Points on Modular Elliptic Curves}, CBMS {\bf 101}, AMS, 2004.

\bibitem{Deb-Kl}
 O. Debarre and  M. Klassen, {\it Points of low degree on smooth plane curves}, J. Reine Angew. Math. {\bf 446} (1994), 81--87. 	 
 
  \bibitem{Del}
 P. Deligne, {\it Formes modulaires et repr\'{e}sentations $\ell$-adiques}, Lecture Notes in Math. 179, Springer-Verlag, 1971, pp.139--172 
 
   \bibitem{SDel}
 P. Deligne and J. P. Serre, { \it Formes modulaires de poids 1},  Ann. Sci. Ec. Norm. Sup. {\bf 7} (1974), 507--530.
 
 
 
 
 
 
 \bibitem{Dec16}
  H. Deconinck, {\it  On the generalized Fermat equation over totally real fields}. Acta Arith. {\bf 173 (3)} (2016), 225--237.

 
 \bibitem{DNS2020}  
 M. Derickx, F. Najman and S. Sikek, {\it Elliptic curve over totally real cubic fields are modular}, Algebra and Number Theory, {\bf 17(7)}, (2020), 1791-1800.
 
Kong, 1993), International Press, 1995, pp. 22--37.


\bibitem{DSbook}
F. Diamond and J. Shurman. A first course in modular forms, volume 228 of Graduate Texts in Mathematics. Springer-Verlag, New York, 2005.

\bibitem{DF2013}
L. Dieulefait, N. Freitas, {\it Fermat-type equations of signature (13,13,p) via Hilbert cuspforms.} Math. Ann. 357 (2013), no. 3, 987--1004.
 

 
 \bibitem{Font-Maz}
J. M.Fontaine and B.Mazur, {\it Geometric Galois representations, in Elliptic Curves, Modular Forms and Fermat’s Last Theorem} (Hong Kong, 1993), International Press, 1995, pp. 41--78.
 
 \bibitem{FKS2020}
 N. Freitas, A. Kraus, and S. Siksek. {\it Class field theory, Diophantine analysis and the asympototic Fermat’s Last Theorem}. Adv. Math. {\bf 363} (2020), 106964.
 
 \bibitem{FHS2015}
  N. Freitas, B. V. Le Hung and S. Siksek, {\it Elliptic curves over real quadratic fields are modular}, Invent. Math. {\bf 201} (2015), no. 1, 159--206.
 
 \bibitem{FS2014} 
  N. Freitas and S. Siksek.{\it  Fermat’s Last Theorem over some small real quadratic fields}, Algebra $\&$ Number Theory{\bf  9} (2014), 875--895.
  
   \bibitem{FS2015} 
  N. Freitas and S. Siksek, {\it The Asymptotic Fermat's Last Theorem for Five-Sixths of Real Quadratic Fields}, Compositio Mathematica, {\bf 151(8)} (2015), 1395--1415.

J. Théor. Nombres Bordeaux {\bf 27} (2015), 67--76.

\bibitem{Fujiwara}
K. Fujiwara. {\it Level optimization in the totally real case} 2006. arXiv: math/0602586  
  

 
 
	Wadsworth and Brooks Cole Advanced Books and Software, Pacific Grove (1990). 
 
 
 	
  	
  	
 
  
 \bibitem{Hida} 
  H. Hida, {\it On abelian varieties with complex multiplication as factors of the Jacobians of
Shimura curves}, Amer. J. Math. {\bf 103 (4)} (1981), 726--776.
 
 
  \bibitem{IKO2020}
 E. I\,{s}ik, Y. Kara, E.  Ozman, {\it  On ternary Diophantine equations of signature
$(p, p, 2)$ over number fields}. Turkish J. Math. {\bf 44}  (2020), no. 4, 1197--1211.
 
 \bibitem{Jarvis}
 F. Jarvis, {\it  Correspondences on Shimura curves and Mazur’s principle at $p$}, Pacific J. Math., {\bf 213 (2)}, 2004, 267--280.
 
 \bibitem{JM2003}
  	F. Jarvis and J. Manoharmayum, {it On the modularity of elliptic curves over totally real fields}, 2003.
 
 \bibitem{JM2008}
 F. Jarvis and J. Manoharmayum, {\it On the modularity of supersingular elliptic curves over certain totally real number fields}, Journal of Number Theory {\bf 128} (2008), no. 3, 589--618.

\bibitem{JM2004}
F. Jarvis and P. Meekin {\it The Fermat equation over $\mathbb{Q}(\sqrt{2})$}, J. Number Theory {\bf 109} (2004), 182--196. 

\bibitem{kal2018} S. Kalyanswamy, {\it Remarks on automorphy of residually dihedral representations}, Math. Res. Lett. (2018), 1285--1304.
 


 

\bibitem{KMBook}
 N. M. Katz and B. Mazur. Arithmetic moduli of elliptic curves, volume 108 of Annals of Mathematics Studies. Princeton University Press, Princeton, NJ, 1985.
 
 \bibitem{MKFJ22}
M. Khawaja  and F. Jarvis,  {\it Fermats ’s Last Theorem over  $\mathbb{Q}(\sqrt{2}, \sqrt{3})$} (2022), https://arxiv.org/abs/2210.03744v2 

 
 
  	
  	
  	 
 
 
 \bibitem{Kraus}
  	A. Kraus, {\it Courbes elliptiques semi-stables et corps quadratiques}, J. Number Theory {\bf 60} (1996), 245--253.
 
 \bibitem{NKSS2022}
 N. Kumar and S. Sahoo, {\it On the solutions of $x^{p}+y^{p} = 2^{r}z^{p}$, $x^{p}+y^{p} =z^{2}$ over totally real fields}, arXiv Pre-print (2022).
 
 
 

\bibitem{Langlands}
R. Langlands, {\it Base Change for $GL_2$}, Ann. of Math. Studies, Princeton University Press 96, 1980.




 
 

 
 \bibitem{PMJ22}
P. Michaud-Jacobs, {\it  Fermat’s Last Theorem and modular curves over real quadratic fields}.
Acta Arith. {\bf 203} (2022), 319--351. 
 
 

 
 
 
 

\bibitem{Raj}
A. Rajaei, {\it On the levels of mod $\ell$ Hilbert modular forms}, J. reine angew. Math. {\bf 537} (2001), 33--65.

 \bibitem{Ribet}
K. A.  Ribet, {\it On modular representations of ${\rm Gal}(\overline{\mathbb{Q}}/\mathbb{Q})$ arising from modular forms}, Invent. Math. {\bf 100} (1990), 431--476.


\bibitem{Ribet97}
 K. A. Ribet, {\it  On the equation  $a^{p}+ 2^{\alpha} b^{p} + c^{p}=0$},  Acta Arith. 
 {\bf 79} (1997), no. 1, 7--1.6

\bibitem{rubin1997}
 K. Rubin, {\it Modularity of mod 5 representations”, pp. 463–474 in Modular forms and Fermat’s last theorem} (Boston, 1995), edited by G. Cornell et al., Springer, 1997
  	
\bibitem{Serre}
J.-P. Serre, {\it Sur les repr\'{e}sentations modulaires de degr\'{e} 2 de ${\rm Gal}(\overline{\mathbb{Q}}/\mathbb{Q})$}, Duke Math J. {\bf 54} (1987), 179--230.  	

\bibitem{ST} 
	T. R. Shemanske, L. H. Walling, {\it Twists of Hilbert modular forms}, Trans. Amer. Math. Soc. {\bf 338 (1)} (1993), 375--403.  

\bibitem{Shi}
	 G. Shimura, {\it The special values of the zeta functions associated with Hilbert modular forms},
	Duke Math. J., {\bf 45(3)} (1978), 637--679. 	

	

 \bibitem{GS71}
G. Shimura, {\rm Introduction to the Arithmetic Theory of Automorphic Functions}, Princeton Univ. Press, Princeton, 1971.

\bibitem{Taylor1989}
R. Taylor, {\it On Galois representations associated to Hilbert modular forms}, Invent. math. {\bf  98} (1989), 265--280.

   
    \bibitem{Taylor-Wiles}
R. Taylor and  A. Wiles, {\it Ring-theoretic properties of certain Hecke algebras}, Ann. Math. {\bf 141} (1995), 553--572.  	
 
 \bibitem{thorne2016} 
 J. Thorne, {\it Automorphy of some residually dihedral Galois representations}, Mathematische Annalen {\bf 364} (2016), no. 1–2, 589--648.
 
 \bibitem{thorne2022} 
 J. Thorne, {\it Elliptic curves over $\mathbb{Q}_{\infty}$ are modular}, to appear in Journal of the European Mathematical Society (2022).
 
 \bibitem{Tunnell}
J.Tunnell, {\it Artin’s conjecture for representations of octahedral type}, Bull.A.M.S.
{\bf 5} (1981), 173--175.
 
  	

 \bibitem{Wiles} 
 A. Wiles, {\it Modular elliptic curves and Fermat’s Last Theorem}, Annals of Mathematics {\bf 141} (1995), no. 3, 443--551. 
 
   
   
   \bibitem{Zhang}
    S.-W. Zhang, {\it  Heights of Heegner points on Shimura curves}, Ann. of Math. {\bf 153} (2001), 27--147.
    
    


	 
 
 


   
 \end{thebibliography}
\end{document}